\documentclass[12pt,a4paper]{amsart}
\usepackage{geometry}
\usepackage{tabls}
\usepackage{url}
\usepackage[utf8]{inputenc}
\usepackage{amsfonts}
\usepackage{amssymb}
\usepackage{mathrsfs}
\usepackage{amsmath}
\usepackage{amsthm}
\usepackage{amscd}
\usepackage{fancyhdr}
\usepackage{enumerate}
\usepackage{paralist}
\usepackage{xypic}
\usepackage{graphicx}
\usepackage{cite}

\usepackage{footmisc}
\usepackage[all,cmtip]{xy}

\numberwithin{equation}{section}

\theoremstyle{plain}

\newtheorem{Def}[equation]{Definition}
\newtheorem{Thm}[equation]{Theorem}
\newtheorem{lem}[equation]{Lemma}

\begin{document}

\title{Iterated integrals, multiple zeta values and Selberg integrals}

\author{Jiangtao Li}

\email{lijiangtao@csu.edu.cn}
\address{Jiangtao Li \\School of Mathematics and Statistics, HNP-LAMA, Central South University, Hunan Province, China}

\begin{abstract}
         Classical multiple zeta values can be viewed as iterated integrals of the differentials $\frac{dt}{t}, \frac{dt}{1-t}$ from $0$ to $1$. In this paper, we reprove Brown's theorem: For $a_i, b_i, c_{ij}\in \mathbb{Z}$, the iterated integral of the form
         \[
         \mathop{\int\cdots \int}\limits_{0<t_1<\cdots<t_N<1}\prod_i t_i^{a_i}(1-t_i)^{b_i} \prod_{i<j}(t_j-t_i)^{c_{ij}}dt_1\cdots dt_N               \]
         is a $\mathbb{Q}$-linear combination of multiple zeta values of weight $\leq N$ if convergent.

         What is more, we show that if $p_i(t), 1\leq i\leq N, $ are in a $\mathbb{Q}\left[t,1/t, 1/(1-t)\right]$-algebra generated by multiple polylogarithms and their dual,
                  and  if $q_{ij}(t), 1\leq i<j\leq N$, are in a $\mathbb{Q}\left[ t,1/t\right]$-algebra generated by logarithm,  then the iterated integral
         \[
          \mathop{\int\cdots \int}\limits_{0<t_1<\cdots<t_N<1}\prod_i p_i(t_i)\prod_{i<j}q_{ij}(t_j-t_i)dt_1\cdots dt_N       \]
          is a $\mathbb{Q}$-linear combination of multiple zeta values.

          As an application of our main results, we show that the coefficients of the Taylor expansions of the Selberg integrals          \[
         \mathop{\int\cdots\int}_{0<t_1<\cdots<t_N<1}f\prod_it_i^{\alpha_i}(1-t_i)^{\beta_i}\prod_{i<j}(t_j-t_i)^{\gamma_{ij}}
           dt_1\cdots dt_N       \]
(with respect to $\alpha_i,\beta_i,\gamma_{ij}$) at the integral points in some product of right half complex plane are $\mathbb{Q}$-linear combinations of multiple zeta values for any $$f\in \mathbb{Q}[t_i, t_i^{-1},(t_i-t_j)^{-1}| 1\leq i\leq N, 1\leq i<j\leq N].$$ This statement generalizes Terasoma's original result.

\end{abstract}

\let\thefootnote\relax\footnotetext{
2020 $\mathnormal{Mathematics} \;\mathnormal{Subject}\;\mathnormal{Classification}$. 11M32, 11G55.\\
$\mathnormal{Keywords:}$  Multiple zeta values, polylogarithms, Selberg integrals. }

\maketitle

\section{Introduction}\label{int}
      Multiple zeta values are defined by
      \[
      \zeta(k_1,\cdots,k_r)=\sum_{0<n_1<\cdots<n_r}\frac{1}{n_1^{k_1}\cdots n_r^{k_r}}, k_1,\cdots, k_{r-1}\geq 1, k_r\geq 2.
      \]
      For multiple zeta value $\zeta(k_1,\cdots,k_r)$, denote by $N=k_1+\cdots+k_r$ and $r$ its weight  and depth respectively. For convenience we view the number $1$ as a multiple zeta value of weight $0$ and depth $0$. Multiple zeta value $\zeta(k_1,\cdots,k_r)$   has an iterated integral representation of the form
      \[
      \zeta(k_1,\cdots,k_r)=\mathop{\int\cdots\int}_{0<t_1<\cdots<t_N<1}\omega_1(t_1)\cdots \omega_N(t_N),
      \]
where $\omega_i(t)=\begin{cases}  \frac{dt}{1-t},&i\in\{1,k_1+1,\cdots,k_1+\cdots+k_{r-1}+1\},\\ \frac{dt}{t},&i\notin\{1,k_1+1,\cdots,k_1+\cdots+k_{r-1}+1\}.\\
\end{cases}
$

Akhilesh \cite{akh} defined double tails of multiple zeta values. They are iterated integrals of the form
\[
\mathop{\int\cdots\int}_{0<t_1<t_2<\cdots<t_N<1}t_1^m(1-t_N)^n\omega_1(t_1)\omega_2(t_2)\cdots \omega_N(t_N),
\]
where $\omega_i\in \Big{\{}\frac{dt}{t},\frac{dt}{1-t}\Big{\}}$. Akhilesh used double tails of multiple zeta values to deduce a very fast algorithm to compute multiple zeta values.

In this paper, we study iterated integrals of the form
  \[
         \mathop{\int\cdots \int}\limits_{0<t_1<\cdots<t_N<1}\prod_i t_i^{a_i}(1-t_i)^{b_i} \prod_{i<j}(t_i-t_j)^{c_{ij}}dt_1\cdots dt_N.                  \]

It is clear that double tails of multiple zeta values are special cases of the above iterated integrals.  Actually Brown \cite{bro}, \cite{bro1} proved the following theorem:
\begin{Thm}\label{brown}(Brown)
For $a_i, b_i, c_{ij}\in \mathbb{Z}$, the iterated integral of the form
         \[
         \mathop{\int\cdots \int}\limits_{0<t_1<\cdots<t_N<1}\prod_i t_i^{a_i}(1-t_i)^{b_i} \prod_{i<j}(t_i-t_j)^{c_{ij}}dt_1\cdots dt_N                 \]
         is a $\mathbb{Q}$-linear combination of multiple zeta values of weight $\leq N$ if convergent.
         \end{Thm}

         In this paper, we will give a different approach to prove the above theorem based on a generalized version of Markarian's main results \cite{mar}. The strategy is: Firstly by changing of variables we view the above integrals as a generalized  version of big zeta value. Secondly, by some trick one can remove the constant terms in the definition of generalized big zeta values. Lastly, by induction on the depth and changing of variables, one reduces the generalized big zeta value to a $\mathbb{Q}$-linear combination of big zeta values.

Ball and Rivoal \cite{BR} introduced the following integral to study irrationality of Riemann zeta values at odd integers:
\[
I_n(z)=\mathop{\int}_{[0,1]^{a+1}}\left( \frac{\Pi_{l=1}^{a+1}x_l^r(1-x_l)}{(z-x_1x_2\cdots x_{a+1})^{2r+1}}\right)^n\frac{dx_1dx_2\cdots dx_{a+1}}{(z-x_1x_2\cdots x_{a+1})^2}, 1\leq r<\frac{a}{2}.
\]
For $z=1$, Ball and Rivoal showed that $I_n(1)$ is a $\mathbb{Q}$-linear combination of Riemann zeta values. By changing of variables $$t_1=x_1x_2\cdots x_{a+1},\cdots,t_i=x_ix_{i+1}\cdots x_{a+1},\cdots t_{a+1}=x_{a+1},$$
it is clear that the integral $I_n(1)$ is a $\mathbb{Q}$-linear combination  of iterated integrals in Theorem \ref{brown}. Dupont\cite{dupont} studied the following integral
\[
\mathop{\int}_{[0,1]^n}\frac{P(x_1,\cdots,x_n)}{(1-x_1\cdots x_n)^N}dx_1\cdots dx_n,
\]
where $P(x_1,\cdots,x_n)$ is a polynomial with rational coefficients. Similarly, one can also show that this kind of integrals  are $\mathbb{Q}$-linear combinations of iterated integrals in Theorem \ref{brown}.

Now we would like to discuss more general kinds of  integrals.

One variable multiple polylogarithms are defined by
\[\mathrm{Li}_{k_1,\cdots,k_r}(t)=\sum_{0<n_1<\cdots<n_r} \frac{t^{n_r}}{n_1^{k_1}\cdots n_r^{k_r}}, |t|<1\]
for $k_1,\cdots,k_r\geq1$. It is clear that $\mathrm{Li}_1(t)=\sum_{n\geq 1}\frac{t^n}{n}=-\mathrm{log}(1-t)$. Let $N=k_1+\cdots+k_r$. For $0<t<1$, it can be written as
 \[
      \mathrm{Li}_{k_1,\cdots,k_r}(t)=\mathop{\int\cdots\int}_{0<s_1<\cdots<s_N<t}\omega_1(s_1)\cdots \omega_N(s_N),
      \]
where $\omega_i(s)=\begin{cases}  \frac{ds}{1-s},&i\in\{1,k_1+1,\cdots,k_1+\cdots+k_{r-1}+1\},\\ \frac{ds}{s},&i\notin\{1,k_1+1,\cdots,k_1+\cdots+k_{r-1}+1\}.\end{cases}$

 Denote by ${\bf P}_{\mathbb{Q}}^{\mathrm{log}}$ the $\mathbb{Q}$-algebra generated by
 \[
 1,\mathrm{Li}_{k_1,\cdots,k_r}(t),\mathrm{Li}_{k_1,\cdots,k_r}(1-t),\forall k_1,\cdots,k_r\geq1.
 \]
 In \cite{li}, the author proved that the iterated integrals
 \[
 \mathop{\int\cdots\int}_{0<t_1<\cdots<t_N<1}f_1(t_1)\omega_1(t_1)\cdots f_N(t_N)\omega_N(t_N), \;f_i(t)\in {\bf P}_{\mathbb{Q}}^{\mathrm{log}}, \omega_i(t)\in \Big{\{}\frac{dt}{t},\frac{dt}{1-t}\Big{\}}.
 \]
 are in the algebra of multiple zeta values if they are convergent (Theorem $1.1$ in \cite{li}).
 Denote by ${\mathbb{Q}}^{\mathrm{log}}[t,\frac{1}{t}]$ the $\mathbb{Q}$-algebra generated by
 \[
 1,\;t,\; \frac{1}{t}, \;\mathrm{Li}_{1}(1-t)=-\mathrm{log}\,t.
 \]
  In this paper, we generalize  Theorem $1.1$ in \cite{li} to the following theorem:
  \begin{Thm}\label{poly}
 For $$p_i(t)\in {\bf P}^{\mathrm{log}}_{\mathbb{Q}}\left[t,\frac{1}{t},\frac{1}{1-t}
 \right], 1\leq i\leq N, $$
 $$
 q_{ij}(t)\in {\mathbb {Q}}^{\mathrm{log}}\left[t,\frac{1}{t}
 \right], 1\leq i<j\leq N,
$$ the iterated integral
 \[
 \mathop{\int\cdots\int}_{0<t_1<\cdots<t_N<1}\prod_i p_i(t_i)\prod_{i<j}q_{ij}(t_j-t_i)dt_1\cdots dt_N,  \]
 is a $\mathbb{Q}$-linear combination of multiple zeta values if it is convergent.
\end{Thm}
        The conditions about the convergence of the above iterated integrals are more complicated than the conditions of Theorem $1.1$ in \cite{li} (see Remark $1.2$ in \cite{li}). We will prove Theorem \ref{poly} by Theorem \ref{brown} and the standard technique in \cite{li}.

        Although Theorem \ref{poly} and Theorem $1.1$ in \cite{li} look very similar, the introduction of the variables $t,\frac{1}{t},\frac{1}{1-t}$ and the terms $q_{ij}(t_i-t_j)$ brings some interesting new phenomenons. Roughly speaking,  iterated integrals of monomials in Theorem $1.1$, \cite{li} are $\mathbb{Q}$-linear combinations of multiple zeta values of pure weight, while iterated integrals of monomials in Theorem \ref{poly} are $\mathbb{Q}$-linear combinations of multiple zeta values of mixed weight.

       For $\alpha_i, \beta_i, \gamma_{ij} \in \mathbb{C}$ and $f\in \mathbb{Q}\left[t_i, \frac{1}{t_i},\frac{1}{t_i-t_j}\right]$, the   Selberg integral is defined as
       \[
       \mathop{\int\cdots\int}_{0<t_1<\cdots<t_N<1}f\prod t_i^{\alpha_i}(1-t_i)^{\beta_i}\prod_{i<j}(t_j-t_i)^{\gamma_{ij}}dt_1\cdots dt_N.
       \]
       By using the theory of Drinfeld associator, Terasoma \cite{tera} proved that: for a suitable choice of $f$, the above Selberg integral is a holomorphic function on $\alpha_i, \beta_i, \gamma_{ij}$ at origin and the coefficients of the Taylor expansion of the above integral with respect to $\alpha_i,\beta_i$ and $\gamma_{ij}$ are $\mathbb{Q}$-linear combinations of multiple zeta values.

       As an application of Theorem \ref{poly}, we prove that
\begin{Thm}\label{beta}
For any $f\in \mathbb{Q}\left[t_i, \frac{1}{t_i},\frac{1}{t_i-t_j}\right]$, the Selberg integral
\[
I(\alpha_i, \beta_i,\gamma_{ij})=\mathop{\int\cdots\int}_{0<t_1<\cdots<t_N<1}f\prod t_i^{\alpha_i}(1-t_i)^{\beta_i}\prod_{i<j}(t_j-t_i)^{\gamma_{ij}}dt_1\cdots dt_N
\]
is a holomorphic function on $$\mathrm{Re}(\alpha_i)>A_i,\mathrm{Re}(\beta_i)>B_i, \mathrm{Re}(\gamma_{ij})>\Gamma_{ij}$$ for some $A_i, B_i, \Gamma_{ij}$ (depends  on $f$). The coefficients of the Taylor expansion of $I(\alpha_i,\beta_i,\gamma_{ij})$ (with respect to $\alpha_i,\beta_i, \gamma_{ij}$) at any integral point in the right half plane $$\mathrm{Re}(\alpha_i)>A_i,\mathrm{Re}(\beta_i)>B_i, \mathrm{Re}(\gamma_{ij})>\Gamma_{ij}$$ are $\mathbb{Q}$-linear combinations of multiple zeta values.
\end{Thm}

Theorem \ref{beta} generalizes Terasoma's original result. Brown (Corollary $8.5$ in \cite{bro}) obtained similar type of results, while the approach here is totally different.

   \section{Generalized big zeta values}\label{mk}

   Markarian \cite{mar} defined the big zeta values and he proved that  the convergent big zeta values are $\mathbb{Q}$-linear combinations of multiple zeta values.

      We review Markarian's theory of big zeta values firstly.
      \begin{Def}
      An element ${\bf e}_{ab}=(r_i)\in \mathbb{Z}^w$ is called a positive $A_w$-root if
      \[
      r_i=\begin{cases}
      1, &a\leq i\leq b\\
      0,& \mathrm{else}.\\
      \end{cases}
      \]
      A $d\times w$-matrix $A=(a_{ij})$ is called basic  if its rank is $d$ and it has no zero columns and all the rows are positive $A_w$-roots.
      \end{Def}
      \begin{Def}
      For a basic matrix $A$, the big zeta value is a formal series
      \[
      Z(A)=\sum_{\substack{n_i\in\mathbb{N}, \\1\leq i\leq d}}\frac{1}{\prod_j(\sum_i a_{ij}n_i)}.
      \]
We call $w,d$ the weight and depth of $Z(A)$ respectively.
      \end{Def}

      In fact, the formal big zeta value $Z(A)$ corresponds to the formal integrals
      \[
      Z(A)=\int_{[0,1]^w}\left(\prod_{j=1}^d\frac{\prod_{ k} x_k^{a_{jk}}}{1-\prod_{ k}x_k^{a_{jk}}}\cdot \prod_{i=1}^w\frac{dx_i}{x_i}\right).
      \]
      Markarian proved the following results.\\
      \begin{Thm} \label{mark}
(Markarian)      Any formal big zeta value is a rational linear combination of formal multiple zeta values of the same weight.
      \end{Thm}
      As a corollary  of Theorem \ref{mark},  we have
      \begin{Thm}\label{spe} (Markarian)
      Any convergent big zeta value is a rational combination of multiple zeta values of the same weight.
      \end{Thm}

      The strategy of Markarian's proof is as follows:\\
 {\bf Step $1$:} For a basic matrix $A$ and a differential operator $D\in \mathbb{Q}\left[ \frac{\partial}{\partial z_1},\cdots, \frac{\partial}{\partial z_d}\right]$ (the trivial operator $1$ is also included here), define
 \[
 Z(A,D)=\sum_{\substack{n_i\in\mathbb{N},\\1\leq i\leq d}}\left(D\frac{1}{\prod_j(\sum_u a_{ij}z_i)}\right)(n_1,\cdots,n_d).
 \]
 By some standard results about hyperplanes (see Proposition 3 in \cite{mar} or Lemma 2 in \cite{BV}), the series $Z(A,D)$ can be written as a $\mathbb{Q}$-linear combination of $Z(S, DD^\prime)$, where $S$ runs over all square submatrices of $A$ of rank $d$.

 {\bf Step $2$:}  Denote by $T_d=\left(a_{ij}\right)_{1\leq i, j\leq d}$ the $d\times d$-matrix satisfying
 \[
 a_{ij}=\begin{cases}
 1, &i\leq j,\\
 0,&i>j.
 \end{cases}
 \]
 By the so-called Orlik-Solomon relations and harmonic product relations, one can reduce $Z(A,D)$ to a $\mathbb{Q}$-linear combination of $Z(T_d, D_{w-d})$ for some homogeneous diffenrential operators $D_{w-d}$. It is clear that $Z(T_d, D_{w-d})$ is a $\mathbb{Q}$-linear formal multiple zeta values of weight $w$. Thus Theorem \ref{mark} is proved.

 Markarian's original results are not enough for our purpose. We need a slight generalization   of Markarian's theorem.
 \begin{Def}
 A $d\times w$-matrix $A=(a_{ij})$ is called quasi-basic  if  it has no zero columns and all the rows are positive $A_w$-roots.
  \end{Def}
  The following lemma is due to Zerbini (Lemma $4.1.2$ in \cite{zer}):
  \begin{lem}\label{ze} (Zerbini)
  For a quasi-basic matrix $A$,  the determinant of its each square sub-matrix  is $0$ or $\pm 1$.
  \end{lem}

  \begin{Def}
  For a quasi-basic matrix $A$ and a row vector $(c_1,\cdots , c_w)\in \mathbb{Z}^w$, the generalized formal big zeta value  is
  \[
  Z(A\,|\,c_1,\cdots,c_w)=\sum_{\substack{n_i\in\mathbb{N}, \\1\leq i\leq d}}\frac{1}{\prod_j(\sum_i a_{ij}n_i+c_j)}.
  \]
  Similarly, we call $w,d$ the weight and depth of $Z(A\,|\, c_1,\cdots, c_w)$ respectively.
  \end{Def}
  \begin{Def}
  For a quasi-basic $d\times w$-matrix $A$ and a differential operator $D\in \mathbb{Q}\left[ \frac{\partial}{\partial z_1},\cdots, \frac{\partial}{\partial z_d}\right]$ (the trivial operator $1$ is also included here), define
 \[
 Z(A,D\,|\, c_1,\cdots,c_w)=\sum_{\substack{n_i\in\mathbb{N},\\1\leq i\leq d}}\left(D\frac{1}{\prod_j(\sum_i a_{ij}z_i+c_j)}\right)(n_1,\cdots,n_d).
 \]
 For $(c_1,\cdots,c_w)=(0,\cdots,0)$, we write
$ Z(A,D)=Z(A,D\,|\,0,\cdots,0)$ for short.    \end{Def}

    For convenience, we review the Orlik-Solomon relations and harmonic product relations for quasi-basic matrices firstly. The proofs are the same as the basic matrices.
 \begin{lem}\label{osr}{(Orlik-Solomon relations)}
 Let $A$ be a quasi-basic matrix and $v_1,\cdots,v_k$  be a set of its distinct columns such that $\sum_i\beta_i v_i=0$ for $0\neq \beta_i\in \mathbb{Q}$. For a differential operator with consatnt coefficients $D$, one has 
 \[
 \sum_{i=1}^k\beta_iZ(\hat{A}_i,D)=0.
 \]
 Here $\hat{A}_i$ is the matrix $A$ with the column $v_i$ excluded.
 \end{lem}
 Denote by 
 \[m_1^*:\mathbb{Q}\left[\frac{\partial}{\partial z_1}, \frac{\partial}{\partial z_2}\right]\rightarrow \mathbb{Q}\left[\frac{\partial}{\partial z_1},  \frac{\partial}{\partial z_2}\right],
 \]   
   \[m_2^*:\mathbb{Q}\left[\frac{\partial}{\partial z_1}, \frac{\partial}{\partial z_2}\right]\rightarrow \mathbb{Q}\left[\frac{\partial}{\partial z_1},  \frac{\partial}{\partial z_2}\right],
 \]  
  \[m_3^*:\mathbb{Q}\left[\frac{\partial}{\partial z_1}, \frac{\partial}{\partial z_2}\right]\rightarrow \mathbb{Q}\left[\frac{\partial}{\partial z}\right],
 \]     
 the linear maps which are induced by  
 \[
 \begin{split}
 &m_1: \mathbb{Q}[z_1,z_2]\ni p(z_1,z_2)\mapsto p(z_1+z_2,z_1)\in \mathbb{Q}[z_1,z_2],\\
&m_2: \mathbb{Q}[z_1,z_2]\ni p(z_1,z_2)\mapsto p(z_1+z_2,z_2)\in \mathbb{Q}[z_1,z_2],\\
 &m_3:\mathbb{Q}\ni p(z)\mapsto p(z_1+z_2)\in\mathbb{Q}[z_1,z_2].\\ 
 \end{split}  \tag{1}\] 
  \begin{lem} (Harmonic product relations)
  Let $A$ be a quasi-basic matrix containing rows ${\bf e}_{ij}$ and ${\bf e}_{(j+1)k}$. Denote by $A^\prime$ the matrix $A$ with these rows excluded.Then we have
  \[
 Z(A,D)= Z([{\bf e}_{ik}; {\bf e}_{ij}; A^\prime], m_1^*D) +Z([{\bf e}_{ik}; {\bf e}_{(j+1)k}; A^\prime], m_2^*D)+Z([{\bf e}_{ik}; A^\prime], m_3^*D), \]
 where $m_i^*$ acts trivially on the subspace generated by rows of $A^\prime$ and acts as in $(1)$ on the subspace generated by the added rows.
   \end{lem}
   
    \begin{Thm}\label{gebig}
    Any  convergent generalized formal big zeta value of weight $w$  and depth $d$ is a $\mathbb{Q}$-linear combination of  multiple zeta values of weight $\leq w$ and depth $\leq d$.
    \end{Thm}

    \noindent{\bf Proof:}
    { \bf Step $1$:} 
 Denote by $d^\prime$ the rank of $A$. By similar standard results about hyperplanes (see Proposition 3 in \cite{mar} or Lemma 2 in \cite{BV}), the series $Z(A,D\,|\, c_1,\cdots, c_w)$ can be written as a $\mathbb{Q}$-linear combination of $Z(S, DD^\prime\,|\, {\bf c^\prime})$ for some $deg \,D^\prime\leq w-d$. Here $S$ runs over all $d\times d^\prime$-submatrices of $A$ of rank $d^\prime$ and  $c^\prime$ is the corresponding $1\times d^{\prime}$-submatrix of $(c_1,\cdots,c_w)$ .
 
  For example,
 \[
 \begin{split}
 &\;\;\;\; \frac{1}{z_1z_2z_3(z_1+z_2+z_3+1)}\\
 &= \frac{1}{z_1z_2z_3}-\frac{1}{z_2z_3(z_1+z_2+z_3+1)} - \frac{1}{z_1z_3(z_1+z_2+z_3+1)} -\frac{1}{z_1z_2(z_1+z_2+z_3+1)}  , \\
 & \;\;\;\; \frac{1}{z_1z_2z_3(z_1+z_2+z_3)}       \\
 &=\frac{1}{z_2z_3(z_1+z_2+z_3)^2}+\frac{1}{z_1z_3(z_1+z_2+z_3)^2}+\frac{1}{z_1z_2(z_1+z_2+z_3)^2} .
  \end{split}
 \]

{\bf Step $2$}:  Denote by $S$ is a $d\times d^\prime$ matrix of rank $d^\prime$, let $S=(s_{ij})_{d\times d^\prime}$. The vector  $(c_1^\prime,\cdots, c_{d^\prime}^\prime)$ satisfies that all the denominators in   $Z(S,D\,|\, c^\prime_1,\cdots, c^\prime_{d^\prime})$ are non-zero. Since $S$ is full rank, by Lemma \ref{ze}, there are a set of  integers $(\mathfrak{T}_1,\cdots, \mathfrak{T}_{d^\prime})$ which satisfy
\[
\sum_{i=1}^{d}s_{ij}\mathfrak{T}_i+c^\prime_j=0,\;j=1,\cdots, d^\prime.
\]
Denote by $(t_1,\cdots,t_{d}$)  a  solution of the above system of linear equations such that $t_i\in \mathbb{Z}, \forall\,i$ and  $\lvert t_1\rvert+\cdots +\lvert t_{d}\rvert$ is minimal among the integer solutions. Then 
 \[
 \begin{split}
 &\;\;\;\;Z(S,D\,|\, c^\prime_1,\cdots,c^\prime_{d^\prime})\\
 &=\sum_{\substack{n_i\in\mathbb{N},\\1\leq i\leq d}}\left(D\frac{1}{\prod_j(\sum_i s_{ij}z_i+c^\prime_j)}\right)(n_1,\cdots,n_{d})\\
 &= \sum_{\substack{n_i\in\mathbb{N},\\1\leq i\leq d}}\left(D\frac{1}{\prod_j(\sum_i s_{ij}(z_i-t_i)+\sum_i s_{ij}t_i+c^\prime_j)}\right)(n_1,\cdots,n_{d})\\
 &= \sum_{\substack{n_i\geq 1-t_i,\\1\leq i\leq d}}\left(D\frac{1}{\prod_j\sum_i s_{ij}z_i}\right)(n_1,\cdots,n_{d}).
  \end{split}
 \]
  Thus one can check that 
  \[
  \begin{split}
  &\;\;\;\;Z(S,D\,|\, c^\prime_1,\cdots,c^\prime_{d^\prime})- \sum_{\substack{n_i\geq 1,\\1\leq i\leq d}}\left(D\frac{1}{\prod_j\sum_i s_{ij}z_i}\right)(n_1,\cdots,n_{d}) 
   \\
  &=\sum_{\substack{n_i\geq 1-t_i,\\1\leq i\leq d}}\left(D\frac{1}{\prod_j\sum_i s_{ij}z_i}\right)(n_1,\cdots,n_{d}) - \sum_{\substack{n_i\geq 1,\\1\leq i\leq d}}\left(D\frac{1}{\prod_j\sum_i s_{ij}z_i}\right)(n_1,\cdots,n_{d}) 
  \end{split}
  \]
  is a $\mathbb{Q}$-linear combination of  generalized big zeta values of  depth $<d$.
  Thus  by induction on the depth $d$, $Z(S,D\,|\, c^\prime_1,\cdots,c^\prime_{d^\prime})$ is reduce to a $\mathbb{Q}$-linear combination of the form $Z(S^\prime, D^{\prime\prime})$. Here $S^\prime$ is a $p\times p^\prime$ quasi-basic matrix of rank $p^\prime$  and $p\leq d$, $\mathrm{deg}\,D^{\prime\prime}+p^\prime\leq \mathrm{deg}\,D+d^\prime$.

 {\bf Step $3$}:  
  We want  to decompose $Z(A,D)$ as a $\mathbb{Q}$-linear combination of
 \[
 Z(U, D^\prime)
 \]
 for some upper triangular block quasi-basic $\delta\times \delta^\prime$-matrix $U$ of rank $\delta^\prime$ (and in every non-zero block all the rows are the same), where $\delta\leq d, \delta^\prime \leq d^{\prime}$. This is achieved by repeated using of harmonic product relations (the same as Procedure $1$ in Section $3$, \cite{mar}). The following lemma is easy to check.
 \begin{lem}\label{most}
 For a quasi-basic $d\times d^\prime$-matrix $A$ of rank $d^{\prime}$, if its $i_1$-th row,  $i_2$-th row, $\cdots$, $i_r$-th row  in $A$ are identical, i.e.
 $$(a_{i_1 1}, a_{i_1 2}, \cdots, a_{i_1 d^{\prime}})=\cdots =(a_{i_r 1}, a_{i_r 2}, \cdots, a_{i_r d^{\prime}})$$ for some $i_1<i_2<\cdots< i_r$,
 then for 
 $$f(z_1,\cdots,z_d)\in \mathbb{Q}[z_1,\cdots, z_d],$$
 there exists some 
 $$ f^\prime( z_1,\cdots, z_{i_2-1}, z_{i_2+1},\cdots,z_{i_r-1},z_{i_r+1},\cdots,  z_{d^\prime})\in \mathbb{Q}[z_i\,|\, i\neq i_2,\cdots, i_r],$$ 
 such that \[
\begin{split}
&\;\;\;\;\;\sum_{n_i\in\mathbb{N}}\frac{f(z_1,\cdots,z_d)}{\prod_j(\sum_i a_{ij}z_i)}(n_1,\cdots, n_d)\\
&=\sum_{\substack{n_i\in \mathbb{N}\\ i\neq i_2,\cdots, i_r}}\binom{z_{i_1}-1}{r-1} \cdot\\
&\frac{f^\prime( z_1,\cdots, z_{i_2-1}, z_{i_2+1},\cdots,z_{i_r-1},z_{i_r+1},\cdots,  z_{d^\prime})}{\prod_j\left(\sum\limits_{i\neq i_2,\cdots,i_r}a_{ij}z_i \right)}\left( n_1,\cdots, n_{i_2-1}, n_{i_2+1},\cdots,n_{i_r-1},n_{i_r+1},\cdots,  n_{d^\prime}\right).
\end{split}
\]
Here $\binom{z_{i_1}-1}{r-1}=\frac{(z_{i_1}-1)\cdots (z_{i_1}-r+1)}{1\cdot 2\cdot \;\cdots\; \cdot (r-1)}$ . \end{lem}

 By using  Lemma \ref{most} to the upper triangular block quasi-basic $\delta \times  \delta^{\prime}$-matrix $U$ of rank $\delta^{\prime}$ for all the identical rows, it follows that $Z(A, D)$ is a $\mathbb{Q}$-linear combination of  series of the form:
 \[
 \sum_{n_1,\cdots,n_{\delta^\prime}\geq 1} \left(z_1^{k_1}\cdots z_{\delta^\prime}^{k_{\delta^\prime}}D^\prime\frac{1}{\prod_{1\leq j\leq \delta^\prime}\left(\sum_{1\leq i\leq \delta^\prime} a^{\prime}_{ij}z_i\right)}\right)\left(n_1,\cdots,n_{\delta^\prime} \right). \tag{2}
 \]
 Here $$k_1,\cdots, k_{\delta^\prime}\geq 0, k_1+\cdots+k_{\delta^\prime}\leq \delta-\delta^\prime,$$
 $$ D^{\prime}\in \mathbb{Q}\left[\frac{\partial}{\partial z_1},\cdots, \frac{\partial}{\partial z_{\delta^\prime}}   \right], \delta^\prime+{deg}\, D^\prime\leq w+{deg}\, D$$ and $A^\prime=\left( a_{ij}^{\prime}\right)$ is an upper triangular {\bf basic} $\delta^\prime \times \delta^\prime$-matrix.

 {\bf Step 4}:
 If $k_1+\cdots+k_{\delta^\prime}=0$ in series $(2)$, we leave the above terms unchanged.
If $k_1+\cdots+k_{\delta^\prime}>0$ in series $(2)$, since ${det}\,(A^\prime)\neq 0$,  we have
 \[
 z_i=\sum_{1\leq j\leq \delta^\prime} p_{ij}\left(\sum_{1\leq i_1\leq \delta^\prime}a^{\prime}_{i_1j}z_{i_1}\right) \tag{3}
 \]
for some $p_{ij}\in\mathbb{Q}$, $1\leq i,j \leq \delta^\prime$.  Thus the series $(2)$ is a $\mathbb{Q}$-linear combination of
\[
\sum_{n_1,\cdots, n_{\delta^\prime}\geq 1}\left(\frac{\prod_{j\in J_2}\left(\sum_{1\leq i\leq \delta^\prime}a_{ij}^{\prime}z_i\right)^{l_j}}{\prod_{j\in J_1}\left(\sum_{1\leq i\leq \delta^\prime}a_{ij}^{\prime}z_i\right)^{l_j}} \right)(n_1,\cdots,n_{\delta^\prime}), \tag{4}
\]
where $$ J_1\cup J_2=\{1, \cdots, \delta^{\prime}\}, J_1\cap J_2=\emptyset,$$
$$\sum_{j\in J_1} l_j\leq \delta^\prime+deg\, D^\prime, \sum_{j\in J_2} l_j<k_1+\cdots+k_{\delta^\prime}.$$
One can apply  Lemma \ref{most}  to the series $(12)$ again. Thus the series $(12)$ is a $\mathbb{Q}$-linear combination of
\[
\sum_{n_1,\cdots, n_s\geq 1}\left(\frac{z_1^{m_1}\cdots z_s^{m_s}}{\prod_{1\leq j\leq s}(\sum_{1\leq i\leq s}b_{ij}z_i)^{p_j}}\right)(n_1,\cdots,n_s), \tag{5}
\]
where $s=\# \,J_1$, $m_1,\cdots,m_s\geq 0, \sum_j p_j=\sum_{j\in J_1}l_j$ and $B=(b_{ij})$ is an upper triangular  {\bf basic} $s\times s$-matrix.
If $m_1+\cdots+m_s>0$ in the series $(5)$, one can repeat the above algorithm. Since $deg\, D+ w$ is a finite number, this algorithm must stop at some step.

{\bf Step $5$}: As a result, the original series $Z(A, D)$    is a $\mathbb{Q}$-linear combination of
 \[
 \sum_{n_1,\cdots,n_{\delta^\prime}\geq 1} \left(D^\prime\frac{1}{\prod_{1\leq j\leq \delta^\prime}\left(\sum_{1\leq i\leq \delta^\prime} a^{\prime}_{ij}z_i\right)}\right)\left(n_1,\cdots,n_{\delta^\prime} \right).
 \]
 Here
 $$ \delta^\prime\leq d, D^{\prime}\in \mathbb{Q}\left[\frac{\partial}{\partial z_1},\cdots, \frac{\partial}{\partial z_{\delta^\prime}}   \right], \delta^\prime+{deg}\, D^\prime\leq w+{deg}\, D$$ and $A^\prime=\left( a_{ij}^{\prime}\right)$ is an upper triangular {\bf basic} $\delta^\prime \times \delta^\prime$-matrix.

 By the  Orlik-Solomon relations, harmonic product relations and  the same way as  Procedure $2$, Section $3$ in \cite{mar}, one can deduce that the series
 \[
 \sum_{n_1,\cdots,n_{\delta^\prime}\geq 1} \left(D^\prime\frac{1}{\prod_{1\leq j\leq \delta^\prime}\left(\sum_{1\leq i\leq \delta^\prime} a^{\prime}_{ij}z_i\right)}\right)\left(n_1,\cdots,n_{\delta^\prime} \right)
  \]
  is actually a $\mathbb{Q}$-linear combination of
  \[
  Z(T_t, E)
  \]
  for some $t\leq \delta^\prime$ and $t+deg\,E\leq \delta^\prime+deg\, D^\prime$. By definition,   $Z(T_t, E)$ is a $\mathbb{Q}$-linear combination of multiple zeta values of weight $\leq t+deg\,E$ and depth $\leq t$.

  In conclusion, any convergent generalized formal big zeta value of weight $w$ and depth $d$ is a $\mathbb{Q}$-linear combination of formal multiple zeta values of weight $\leq w$ and depth $\leq d$. $\hfill\Box$\\

\section{Proof of main theorem and application}

Brown \cite{bro1}  sketched a proof of the following theorem:
\begin{Thm}\label{or}(Brown) For $a_i, b_i,c_{ij}\in\mathbb{Z}$ and $l\geq 1$, the iterated integral
\[
\mathop{\int\cdots\int}_{0<t_1<\cdots<t_l<1}\prod_i t_i^{a_i}(1-t_i)^{b_i}\prod_{i<j}(t_i-t_j)^{c_{ij}}dt_1\cdots dt_l
\]
is a $\mathbb{Q}$-linear combination of multiple zeta values of weight $\leq l$.
\end{Thm}

Now we  give the detailed proof of Brown's theorem based on the results in Section \ref{mk}.\\
{\bf Proof of Theorem \ref{or}:}
The original iterated integral is a $\mathbb{Q}$-linear combination of the following iterated integrals
\[
\mathop{\int\cdots\int}_{0<t_1<\cdots<t_l<1}\prod_i t_i^{a^\prime_i}(1-t_i)^{b^\prime_i}\prod_{i<j}(t_i-t_j)^{c^\prime_{ij}}dt_1\cdots dt_l,
\]
where $b_i^\prime, c_{ij}^\prime\leq 0$.
Thus it suffices to prove Theorem \ref{or} in the cases $b^\prime_i,c^\prime_{ij}\leq 0$.

Now we assume  $b_i,c_{ij}\leq 0$. Let $t_1=x_1x_2\cdots x_l, t_2=x_2\cdots x_l, \cdots , t_l=x_l$, then
\[
dt_1dt_2\cdots dt_l=x_2x_3^2\cdots x_l^{l-1}dx_1dx_2\cdots dx_l,
\]
\[
\begin{split}
&\;\;\mathop{\int\cdots\int}_{0<t_1<\cdots<t_l<1}\prod_i t_i^{a_i}(1-t_i)^{b_i}\prod_{i<j}(t_i-t_j)^{c_{ij}}dt_1\cdots dt_l
\\
&=\mathop{\int}_{[0,1]^l} \prod_ix_i^{i-1+\sum_{n\leq i}a_n}\prod_{i<j}(x_j\cdots x_l)^{c_{ij}} \prod_i(1-x_i\cdots x_l)^{b_i}\prod_{i<j}(x_i\cdots x_{j-1}-1)^{c_{ij}}dx_1\cdots dx_l.
\end{split}
\]
As a result, the iterated integral
\[
\mathop{\int\cdots\int}_{0<t_1<\cdots<t_l<1}\prod_i t_i^{a_i}(1-t_i)^{b_i}\prod_{i<j}(t_i-t_j)^{c_{ij}}dt_1\cdots dt_l
\]
is equal to
\[
C\cdot\int_{[0,1]^l}\prod_i\left(\frac{1-x_i\cdots x_l}{x_i\cdots x_l}\right)^{b_i}\prod_{i<j}\left(\frac{1-x_i\cdots x_{j-1}}{x_i\cdots x_{j-1}}\right)^{c_{ij}}\prod_i x_i^{f_i}\cdot \frac{dx_i}{x_i} \tag{6}
\]
for some $C\in \{\pm 1\}$ and $f_i\in \mathbb{Z}$.

We use the notations in Section \ref{mk}. Let $$d=\sum_i \lvert b_i\rvert +\sum_{i<j}\lvert c_{ij}\rvert, \,w=l.$$
The $d$ positive $A_w$-roots
\[
\underbrace{{\bf e}_{iw}, \cdots , {\bf e}_{iw}}_{\mid b_i\mid},\;  1\leq i\leq l,
\]
\[
\underbrace{{\bf e}_{i j-1},\cdots, {\bf e}_{i j-1}}_{\mid c_{ij}\mid}, \;1\leq i<j\leq l
\]
form a $d\times w$ matrix $A$.
By using the formula
$\frac{x}{1-x}=\sum_{n\geq 1} x^n$,
it is easy to check that
\[
\int_{[0,1]^l}\prod_i\left(\frac{1-x_i\cdots x_l}{x_i\cdots x_l}\right)^{b_i}\prod_{i<j}\left(\frac{1-x_i\cdots x_{j-1}}{x_i\cdots x_{j-1}}\right)^{c_{ij}}\prod_i x_i^{f_i}\cdot \frac{dx_i}{x_i}=Z(A\,|\, f_1,\cdots, f_w).
\]
If there are no non-zero columns in the matrix $A$, then $A$ is quasi-basic. If there are exactly  $b$ zero columns in $A$, then it is easy to check that these zero columns contribute an independent factor of 
the form 
\[
\mathop{\int}_{[0,1]^b}x_{i_1}^{f_{i_1-1}}\cdots x_{i_b}^{f_{i_b-1}}dx_{i_1}\cdots dx_{i_b}.
\]  
As a result, 
\[
Z(A\,|\, f_1,\cdots, f_w)=\mathop{\int}_{[0,1]^b}x_{i_1}^{f_{i_1-1}}\cdots x_{i_b}^{f_{i_b-1}}dx_{i_1}\cdots dx_{i_b}\cdot Z(B\,|\, p_1,\cdots, p_{w^\prime})
\]
for some quasi-basic $d\times w^\prime$ matrix $B$, where $w^\prime=w-b$.

 From Theorem \ref{gebig},  $Z(A\, |\, f_1,\cdots, f_w)$  is a $\mathbb{Q}$-linear combination of formal  multiple zeta values of weight $\leq l$ and depth $\leq d$.
Thus if the original iterated integral is convergent, it is a $\mathbb{Q}$-linear combination of  multiple zeta values of weight $\leq N$.
$\hfill\Box$\\

\section{Multiple polylogarithms and Selberg integrals}
     In this section, we prove Theorem \ref{poly} by using Theorem \ref{or} and the skills in \cite{li}. As an application of Theorem \ref{poly},  we show that the coefficients of the Taylor expansion of Selberg integrals at integral points are in the algebra of multiple zeta values.

     For $k_1,\cdots, k_r\geq 1$, the multiple polylogarithm $\mathrm{Li}_{k_1,\cdots, k_r}(t)$ is defined by $$\mathrm{Li}_{k_1,\cdots, k_r}(t)=\sum_{0<n_1<\cdots<n_r}\frac{t^{n_r}}{n_1^{k_1}\cdots n_2^{k_r}}, \lvert t\rvert<1.$$
         If $k_r\geq 2$, the evaluation of $\mathrm{Li}_{k_1,\cdots, k_r}(t)$   at $t=1$ is $\zeta(k_1,\cdots,k_r)$.

    Denote by $S_{m+n}$ the group of permutations on the set $\{1,2,\cdots,m+n\}$. From the theory of iterated integrals \cite{chen}, the following lemma is true:
     \begin{lem}\label{ite}(Chen) (i) For any $\omega_i$, $i=1,\cdots,m+n$,
    \[
    \begin{split}
    &\;\;\;\;\mathop{\int\cdots\int}_{a<t_1<\cdots<t_m<b}\omega_1(t_1)\cdots \omega_{m}(t_m)\cdot\mathop{\int\cdots\int}_{a<t_{m+1}<\cdots<t_{m+n}<b}\omega_{m+1}(t_{m+1})\cdots\omega_{m+n}(t_{m+n})\\
    &=\sum_{\sigma\in {Sh}(m,n)}\mathop{\int\cdots\int}_{a<t_{\sigma(1)}<t_{\sigma(2)}\cdots<t_{\sigma(m+n)}<b}\omega_1(t_1)\omega_2(t_2)\cdots\omega_{m+n}(t_{m+n}),  \\
    \end{split}
    \]
    where elements of ${Sh}(m,n)\subseteq S_{m+n}$ satisfy
    \[
    \sigma^{-1}(1)<\sigma^{-1}(2)<\cdots<\sigma^{-1}(m)
    \]
    and
    \[
    \sigma^{-1}(m+1)<\sigma^{-1}(m+2)<\cdots< \sigma^{-1}(m+n).
    \]
    (ii) For $a<c<b$,
    \[
    \begin{split}
    &\;\;\;\;\mathop{\int\cdots\int}_{a<t_1<\cdots<t_n<b}\omega_1(t_1)\cdots \omega_{n}(t_n)\\
    &=\mathop{\int\cdots\int}_{a<t_1<\cdots<t_n<c}\omega_1(t_1)\cdots \omega_{n}(t_n)+\mathop{\int\cdots\int}_{c<t_1<\cdots<t_n<b}\omega_1(t_1)\cdots \omega_{n}(t_n),  \\
    &+\sum_{1\leq i\leq n-1}\mathop{\int\cdots\int}_{a<t_{1}<\cdots<t_{i}<c}\omega_1(t_1)\cdots\omega_{i}(t_{i})\cdot \mathop{\int\cdots\int}_{c<t_{i+1}<\cdots<t_{n}<b}\omega_1(t_{i+1})\cdots\omega_{n}(t_{n})\\
       \end{split}
    \]

     \end{lem}

     {\bf Proof of Theorem \ref{poly}:} For $$p_i(t) \in {\bf P}^{\mathrm{log}}_{\mathbb{Q}}\left[t,\frac{1}{t},\frac{1}{1-t}\right], 1\leq i\leq N, q_{ij}(t)\in {\bf Q}^{\mathrm{log}}_{\mathbb{Q}}\left[t,\frac{1}{t}\right],1\leq i<j\leq N.$$
It suffices to show that if $p_1(t),\cdots, p_N(t)$ are monomials in ${\bf P}_{\mathbb{Q}}^{\mathrm{log}}[t,\frac{1}{t},\frac{1}{1-t}]$ and $q_{ij}(t), 1\leq i<j\leq N$,  are monomials in ${\mathbb Q}^{\mathrm{log}}[t,\frac{1}{t}]$, the iterated integral
     \[
     \mathop{\int\cdots\int}_{0<t_1<\cdots<t_N<1}\prod_i p_i(t_i)\prod_{i<j}q_{ij}(t_j-t_i)dt_1\cdots dt_N     \]
     is a $\mathbb{Q}$-linear combination of multiple zeta values. In this case, $p_i(t)$ can be written as the form
     \[
     \frac{1}{t^{m_i} {(1-t)}^{n_i}}\prod_{u\geq 1}\prod_{k_1,\cdots,k_u\geq 1}\left(\mathrm{Li}_{k_1,\cdots,k_u}(t)\right)^{p_{k_1,\cdots,k_u}} \cdot \left(\mathrm{Li}_{k_1,\cdots,k_u}(1-t)\right)^{q_{k_1,\cdots,k_u}}, \tag{7}
     \]
     where $m_i,n_i\in \mathbb{Z}$, $p_{k_1,\cdots,k_u}, q_{k_1,\cdots,k_u}\geq 0$.
     Similarly $q_{ij}(t)$ can be written as the form
      \[
     \frac{1}{t^{m_{ij}}}\left( \mathrm{Li}_1(1-t)\right)^{l_{ij}},   \tag{8}  \]
     where $m_{ij}\in\mathbb{Z}$, $ l_{ij}\geq 0$.
     For $0<t<1$,
     \[
     \mathrm{Li}_{k_1,\cdots,k_u}(t)=\mathop{\int\cdots\int}_{0<t_1<\cdots<t_{k_1+\cdots+k_u}<t}\omega_1(t_1)\cdots\omega_{k_1+\cdots+k_u}(t_{k_1+\cdots+k_u})
     \]
     and
       \[
       \begin{split}
     &\;\;\;\;\;\mathrm{Li}_{k_1,\cdots,k_u}(1-t)\\
     &=\mathop{\int\cdots\int}_{0<t_1<\cdots<t_{k_1+\cdots+k_u}<1-t}\omega_1(t_1)\cdots\omega_{k_1+\cdots+k_u}(t_{k_1+\cdots+k_u})\\
     &=\mathop{\int\cdots\int}_{t<t_{k_1+\cdots+k_u}<\cdots<t_1<1}\omega_1(1-t_1)\cdots\omega_{k_1+\cdots+k_u}(1-t_{k_1+\cdots+k_u}) \\
     &=\mathop{\int\cdots\int}_{t<t_{k_1+\cdots+k_u}<\cdots<t_1<1}\omega_1^{\prime}(t_1)\cdots\omega_{k_1+\cdots+k_u}^{\prime}(t_{k_1+\cdots+k_u}) \\
          \end{split}
         \]
         for some $\omega_s(t),\omega^{\prime}_s(t)\in \big{\{}\frac{dt}{t},\frac{dt}{1-t}\big{\}}$, $ s=1,\cdots, k_1+\cdots+k_u$.

By Lemma \ref{ite}, $(i)$, monomials of type $(7)$ can be written as $\mathbb{Q}$-linear combinations of the following integrals
\[
\frac{1}{t^{m_i}(1-t)^{n_i}}\mathop{\int\cdots\int}_{0<u_1<\cdots<u_A<t}\omega_1(u_1)\cdots\omega_A(u_A)\cdot\mathop{\int\cdots\int}_{t<u_{A+1}<\cdots<u_{A+B}<1}\omega_{A+1}(u_{A+1})\cdots \omega_{A+B}(u_{A+B}),
\]
where $A, B\geq 1, \omega_i(u)\in \big{\{}\frac{du}{u},\frac{du}{1-u}\big{\}}, i=1,\cdots, A+B$ and $m_i, n_i\in\mathbb{Z}$.

By Lemma \ref{ite}, $(i)$, monomials of type $(8)$ can be written as $\mathbb{Q}$-linear combinations of the following integrals
\[
\begin{split}
&\frac{1}{(t_j-t_i)^{m_{ij}}}\mathop{\int\cdots \int}_{t_j-t_i<u_1<\cdots<u_{\Gamma}<1}\omega_1(u_1)\cdots \omega_{\Gamma}(u_{\Gamma}),\\
\end{split}
\]
where $\Gamma\geq 1$, $\omega_i(u)=\frac{du}{u},  \forall\, i=1,\cdots, \Gamma$ and $m_{ij}\in \mathbb{Z}$.
By Lemma \ref{ite}, $(ii)$, we have
\[
\begin{split}
&\;\;\;\;\mathop{\int\cdots \int}_{t_j-t_i<u_1<\cdots<u_{\Gamma}<1}\frac{du_1}{u_1}\cdots \frac{du_{\Gamma}}{u_{\Gamma}}\\
&=\mathop{\int\cdots \int}_{t_j-t_i<u_1<\cdots<u_{\Gamma}<1-t_i}\frac{du_1}{u_1}\cdots \frac{du_{\Gamma}}{u_{\Gamma}}+\mathop{\int\cdots\int}_{1-t_i<u_1<\cdots<u_{\Gamma}<1}\frac{du_1}{u_1}\cdots \frac{du_{\Gamma}}{u_{\Gamma}} \\
&+\sum_{1\leq l\leq \Gamma-1}\mathop{\int\cdots\int}_{t_j-t_i<u_{1}<\cdots<u_{l}<1-t_i}\frac{du_1}{u_1}\cdots \frac{du_{l}}{u_{l}}\cdot \mathop{\int\cdots\int}_{1-t_i<u_{l+1}<\cdots<u_{\Gamma}<1} \frac{du_{l+1}}{u_{l+1}}\cdots \frac{du_{\Gamma}}{u_{\Gamma} }    \\
    &=\mathop{\int\cdots \int}_{t_j<u_1<\cdots<u_{\Gamma}<1}\frac{du_1}{u_1-t_i}\cdots \frac{du_{\Gamma}}{u_{\Gamma}-t_i} +\mathop{\int\cdots\int}_{0<u_1<\cdots<u_{\Gamma}<t_i}\frac{du_1}{1-u_1}\cdots \frac{du_{\Gamma}}{1-u_{\Gamma}}   \\
&+\sum_{1\leq l\leq \Gamma-1}\mathop{\int\cdots\int}_{t_j<u_{1}<\cdots<u_{l}<1}\frac{du_1}{u_1-t_i}\cdots \frac{du_{l}}{u_{l}-t_i}\cdot \mathop{\int\cdots\int}_{0<u_{l+1}<\cdots<u_{\Gamma}<t_i} \frac{du_{l+1}}{1-u_{l+1}}\cdots \frac{du_{\Gamma}}{1-u_{\Gamma} }    \\
   \end{split}
\]

As a result, the iterated integral
      \[
     \mathop{\int\cdots\int}_{0<t_1<\cdots<t_N<1}\prod_i p_i(t_i)\prod_{i<j}q_{ij}(t_j-t_i)dt_1\cdots dt_N     \]
     is a  $\mathbb{Q}$-linear combination of the following integrals
     \[
     \mathop{\int}_{\Sigma}\prod_{i=1}^Ct_i^{m_i}(1-t_i)^{n_i}\prod_{1\leq i<j\leq C}(t_j-t_i)^{l_{ij}}dt_1\cdots dt_C.
          \]
     Here  $C>0, m_i, n_i,l_{ij}\in \mathbb{Z}$ for  $i=1,\cdots,C$, $\Sigma\subseteq [0,1]^C$ and $\Sigma$ is defined by  conditions of the form $t_i<t_j$ for some $i,j$ (usually many pairs of  $(i,j)$).
     Since the set of real numbers is a total ordered set, $\Sigma$ can be written as a disjoint union of the set
     \[
     \big{\{}(t_1,\cdots,t_N)\mid 0<t_{\sigma(1)}<\cdots<t_{\sigma(C)}<1\big{\}}
     \]
     for some $\sigma\in S_C$ (plus some lower dimension terms).   Thus Theorem \ref{poly} follows from Theorem \ref{brown}.      $\hfill\Box$\\

Now we are ready to prove Theorem \ref{beta}.

It is clear that for the Selberg integral
\[
I(\alpha_i, \beta_i,\gamma_{ij})=\mathop{\int\cdots\int}_{0<t_1<\cdots<t_N<1}f\prod t_i^{\alpha_i}(1-t_i)^{\beta_i}\prod_{i<j}(t_j-t_i)^{\gamma_{ij}}dt_1\cdots dt_N,
\]
the first statement is true. From the Taylor expansion of $I(\alpha_i,\beta_i, \gamma_{ij})$ at the integral points, we have
\[
\begin{split}
&\;\;\;\;I(\alpha_i, \beta_i,\gamma_{ij})
\\&=\sum_{a_i, b_i, c_{ij}\geq 0} \left[\prod_i \left(\frac{\partial}{\partial \alpha_i}\right)^{a_i}\left(\frac{\partial}{\partial \beta_i}\right)^{b_i}
\prod_{i<j}\left(\frac{\partial}{\partial \gamma_{ij}}\right)^{c_{ij}}
  I(\alpha_i, \beta_i,\gamma_{ij}) \right]_{\substack{(\alpha_i, \beta_i)=(A_i, B_i)\\(\gamma_{ij})=(C_{ij})}}\\
&\;\;\;\;\;\;\;\;\;\;\;\;\;\;\;\;\;\;\;\;\;\cdot\prod_i\frac{(\alpha_i-A_i)^{a_i}}{a_i!}\frac{(\beta_i-B_i)^{b_i}}{b_i!}\prod_{i<j} \frac{(\gamma_{ij}-C_{ij})^{c_{ij}}}{c_{ij}!}\\
&=\sum_{a_i,b_i,c_{ij}\geq 0}\\
&\mathop{\int\cdots\int}_{0<t_1<\cdots<t_N<1} f\prod_i t_i^{A_i}(1-t_i)^{B_i}\mathrm{log}^{a_i}t_i\mathrm{log}^{b_i}(1-t_i)\prod_{i<j}(t_i-t_j)^{C_{ij}} \mathrm{log}^{c_{ij}}(t_i-t_j)   dt_1\cdots dt_N\\
&\;\;\;\;\;\;\;\;\;\;\;\;\;\;\;\;\;\;\;\;\;\cdot\prod_i\frac{(\alpha_i-A_i)^{a_i}}{a_i!}\frac{(\beta_i-B_i)^{b_i}}{b_i!}\prod_{i<j} \frac{(\gamma_{ij}-C_{ij})^{c_{ij}}}{c_{ij}!}.\\
\end{split}
\]
Thus Theorem \ref{beta} is equivalent to that for some fixed $A_i, B_i , C_{ij}\in\mathbb{Z}$ and for every $a_i, b_i, c_{ij}\geq 0$, the iterated integral
\[
\mathop{\int\cdots\int}_{0<t_1<\cdots<t_N<1} f\prod_i t_i^{A_i}(1-t_i)^{B_i}\mathrm{log}^{a_i}t_i\mathrm{log}^{b_i}(1-t_i)\prod_{i<j}(t_i-t_j)^{C_{ij}} \mathrm{log}^{c_{ij}}(t_i-t_j)   dt_1\cdots dt_N
\]
is a $\mathbb{Q}$-linear combination of multiple zeta values.
 Since $\mathrm{Li}_1(1-t)=-\mathrm{log}(t), t\in [0,1]$, the above statement follows immediately from Theorem \ref{poly}.

\section*{Acknowledgements}
       The author wants to thank the referee for the detailed comments to improve this paper.  The author is supported by the National Natural Science Foundation of China (Grant No. 12201642).

\end{document}